\documentclass{article} 
\usepackage{graphicx} 
\begin{document}
\setcounter{page}{1} 
\begin{center}
{\Large\bf An exact closed walks formula for the 
complexity of regular graphs and some related bounds } 
\end{center}
\medbreak
\centerline{by}
\medbreak
\centerline{Gregory P. Constantine}

\centerline{School of Computer Science}

\centerline{Georgia Institute of Technology}

\centerline{Atlanta, GA 30332 }

\medbreak
\centerline{and}
\medbreak
\centerline{Gregory C. Magda}

\centerline{Department of Mathematics}

\centerline{University of Pittsburgh}

\centerline{Pittsburgh, PA 15260}
\vskip1cm 
\centerline{\bf ABSTRACT}
 
\medbreak\noindent The complexity of a graph is the 
number of its labeled spanning trees.  In this work 
complexity is studied in settings that admit regular 
graphs.  An exact formula is established linking 
complexity of the complement of a regular graph to 
numbers of closed walks in the graph by way of an 
alternating series.  Some consequences of this 
result yield infinite classes of lower and upper bounds 
on the complexity of such graphs.  Applications of these 
mathematical results to biological problems on neuronal 
activity are described.  

\noindent {\em Proposed }
{\em running head:\/}  Complexity and closed walks 
\vskip1cm 
\footnoterule\noindent Funded under NIH grant 
RO1-HL-076157 and NSF award ID 2424684 

\noindent gmc@pitt.edu 
  
\newpage
\section*{\Large\bf1.  Motivation and preliminaries} 
\medbreak\noindent The motivating biological problem is to 
turn on all neurons in a brain, or part of a brain, by 
starting with a small subset of active neurons.  We 
view this activity as having {\em local component\/}s, which we 
want to turn on as fast as possible, and {\em global links }
between the local components which serve the purpose 
of efficiently integrating the local components in such a 
way that the entire brain becomes active as quickly as 
possible.  Further detailed information is found in $[14]$ 
and [13].  Intuitively we thus seek to determine those 
neuronal configurations, viewed as abstract networks, 
that spread the information most efficiently (fastest 
possible) to the whole brain.  We focuss initially on 
modeling the local components and start by making some 
simplifying assumptions.  The basic working hypothesis 
is that a neuron is activated by receiving input from at 
least $t$ already active neurons connected to it.  Initially 
we make the assumption that the underlying graph that 
connects the neurons is regular.  Since we want a quick 
spread to activate the whole local area, it is intuitive 
that the best way of doing this is to avoid having short 
cycles, like triangles or squares, in the regular graph.  
If we have $n$ neurons, each of degree $d,$ the emerging 
optimization strategy is that we want to first restrict 
to having a minimal number of triangles, then among 
this subset of regular graphs to seek those that have a 
minimum number of closed walks of length 4 (like 
4-cycles), and proceed sequencially to closed walks of 
higher order.  The point of this paper is to establish a 
mathematical connection between the choice strategy we 
just described and regular graphs of degree $d$ that have 
a maximum number of spanning trees.  We describe next, 
in some detail, measures of the spread of neuronal 
activity.  \medbreak\noindent Imagine for a moment that 
the vertices of the graph (or digraph) are neurons and 
any existing edge transmits information form one neuron 
to another.  We start with a set $S$ of neurons, which 
we call $active,$ and a startup treshold $t,$ which is a 
natural number.  The spreading of neuronal activity is 
described next.  This is subject to some restrictions 
formulated in terms of {\em Steps}, which we now describe.  
\vskip.3cm\noindent
$Step$ 0:  Start with a set $S=S_0$ of vertices of the 
digraph $G$ and a natural number $t.$ We call elements of 
$S$ active vertices.  [Imagine that you hold the active 
vertices in your left hand, and the other vertices in 
your right hand.]  Color any edge emanating from $S$ red.  
\vskip.3cm\noindent
$Step$ 1:  Acquire vertex $v$, held in your right hand, if $v$ 
has $t$ or more red arrows pointing to it.  Move all 
acquired vertices to your left hand.  Call the set of 
vertices you now hold in your left hand $S_1.$ Color all 
edges emanating from $S_1$ red.  
\vskip.3cm\noindent
The general step is as follows.  We are in posession of 
$S_{i-1}$ with all edges emanating from it colored red.  
\vskip.3cm\noindent
$Step$ $i:$ Acquire vertex $v,$ held in your right hand, if it 
has $t$ or more red arrows pointing to it.  Move all 
acquired vertices to your left hand.  Call the set of 
vertices you now hold in your left hand $S_i.$ Color all 
edges emanating from $S_i$ red.  
\vskip.3cm\noindent
Evidently $S=S_0\subseteq S_1\subseteq\cdots\subseteq S_i\subseteq
\cdots .$ As we keep 
increasing $i$, the following will (obviously) always occur:  
the number of vertices in your right hand becomes 
stationary; that is, $\exists m$ such that, at Step $i,$ for all $
i\geq m$ 
the number of vertices in your right hand remains 
constant.  If your right hand becomes empty for a 
sufficiently large $i$ we say that the network is in 
{\em synchrony}.  [You are now holding the whole network in 
your left hand -- hence all vertices of the network 
became active.]  We denote by $i^{*}[=i^{*}(S,t)]$ the smallest 
$i$ such that at Step $i$ the network is in synchrony.  
Typically a network cannot be brought to synchrony 
(starting with an incipient set $S=S_0$ and $t)$, and we 
convene to write $i^{*}=i^{*}(S,t)=\infty$ in such a case.  
\medbreak\noindent Fix $t.$ In a graph (or digraph) $G,$ let 
$S=S_0$ be a set with $k$ vertices, which we call a 
$k-$subset.  Write $i^{*}(S)$ for $i^{*}(S,t).$ We introduce the 
following measures of synchrony for $G$ and $k.$ 
\vskip.3cm\noindent
The ratio $p_k(G)=$(number of $k-$subsets $S$ that bring $G$ to 
synchrony)/(number of all $k-$subsets) signifies the 
probability of bringing digraph $G$ to synchrony from a 
randomly chosen $k-$subset.  Generally we are interested 
in identifying digraphs with large $p_k$.  It might also be 
observed that there are many instances when a digraph 
has a large $p_k$ but the number of steps required to 
obtain synchrony are generally quite large, which is not 
so good.  We could tune this up by defining another 
measure $e_k(G)$, which we call $synchrony$ $efficiency,$ as 
follows:  
\vskip.3cm 
${n\choose k}e_k(G)=\sum_S(i^{*}(S))^{-1}.$ 
\vskip.3cm\noindent
Observe that when $S$ does not induce synchrony, 
$i^{*}(S)=\infty ,$ and we simply add a zero to the sum.  
Intuitively, efficiency $e_k$ yields the average speed to the 
synchrony of $G$ across all $k-$subsets.  High values of $e_k$ 
are typically good, since the synchrony is then speedily 
restored.  We did not see the concept of synchrony 
efficiency used in the network optimization literature so 
far.  Graph theoretic preliminaries are introduced next.  
\medbreak\noindent The graphs we work with are finite, 
loopless, undirected, and without multiple edges.  By the 
{\em order\/} of a graph we understand the number of its 
vertices, and by {\em size\/} the number of its edges.  A graph 
is called $regular$ if the degrees of its vertices are equal.  
Standard terminology is used and we assume that the 
reader is familiar with such notions as path, graph 
connectivity, tree and spanning tree, adjacency matrix 
and the Laplacian; see [6, 10].  For clarity we also 
remind that a $walk$ of length $k$ (or $k-$walk) is a sequence 
of vertices and edges $v_1e_1v_2e_2\cdots v_ke_kv_{k+1}$, where $
e_j$ is the 
edge joining vertices $v_j$ and $v_{j+1}.$ Vertices and edges 
may be repeated in this sequence.  The walk is {\em closed\/} if 
$v_1=v_{k+1}.$ A {\em m-cycle\/} is a sequence of vertices and edges 
$v_1e_1v_2e_2\cdots v_me_mv_{m+1},$ where all vertices $v_i$ are distinct 
except for $v_1$ and $v_{m+1}$ which are the same; $m\geq 2.$ A 
$triangle$ is a $3-$cycle; it is also a closed 3-walk.  
\vskip.5cm\noindent
Denote by $D$ the diagonal matrix with the degrees of the 
vertices of graph $G$ as entries (written always in the 
same fixed order), by $A$ the adjacency matrix and by 
$L=D-A$ the Laplacian.  The {\em complexity of graph} $G$ is 
the number of (labeled) spanning trees of $G;$ we denote 
it by $t(G).$ The {\em log-complexity of} $G$ is defined as 
$log(t(G)).$ We remind the reader of a few 
well-known results, see [6, 1] and [10], that we rely 
on and use freely in this article:  \medbreak\noindent{\bf 1.} The $
(i,j)^{th}$ entry of $A^r$ is equal to the number of walks 
with $r$ edges staring at vertex $v_i$ and ending at vertex 
$v_j$.  In particular, the number of closed $r-$walks at 
vertex $v_i$ is the $(i,i)^{th}$ entry of $A^r.$ Consequently $tr(A^r)$, 
the trace of $A^r$, is equal to $w_r(G)$, the total number of closed 
$r-$walks in graph $G$.  \medbreak\noindent{\bf 2.} If $\lambda_1\geq\lambda_
2\geq\ldots\geq\lambda_n$ are 
the eigenvalues of the $n\times n$ adjacency matrix $A,$ then 
$tr(A^r)=\sum_{i=1}^n\lambda_i^r.$ \medbreak\noindent{\bf 3.}  If the graph is 
regular of degree $d,$ then $L=dI-A,$ with $I$ denoting the 
identity matrix.  Furthermore, the eigenvalues $\mu_i$ of $L$ 
may be written in this case as $\mu_i=d-\lambda_i,$ $1\leq i\leq 
n.$ 
Since the row sums of $L$ are always 0, we have $\mu_1=0.$ 
\medbreak\noindent{\bf 4.} It is an easy consequence of Kirchhoff's 
theorem that {\em the complexity of graph $G$ is equal to $\frac 
1n\mu_2\mu_3\cdots\mu_n$ where $n$ 
is the $order$ of $G$ and $(0=)\mu_1\leq\mu_2\leq
\mu_3\leq\ldots\leq\mu_n$ are the 
eigenvalues of the Laplacian $L$ of $G.$} \medbreak\noindent
The complement $\bar {G}$ of graph 
$G$ is the graph in which $e$ is an edge if $e$ is not an edge 
in $G.$ We denote by $\bar {A}$ and $\bar {L}$ the adjacency matrix and 
the Laplacian of $\bar {G}.$ Let $I$ be the identity matrix and $
J$ 
be the square matrix with all entries equal to 1.  
Evidently $A+\bar {A}=J-I$ and $L+\bar {L}=nI-J.$ These equalities 
allow us to immediately conclude as follows:  
\medbreak\noindent{\bf 5.} The eigenvalues of $\bar {L}$ are $\bar{
\mu}_i=n-\mu_i,$ 
$2\leq i\leq n$ and $\bar{\mu}_1=0.$ In view of {\bf 4.} we have 
$t(\bar {G})=\frac 1n\prod_{i=2}^n\bar{\mu}_i=n^{n-2}\prod_{i=2}^
n(1-\frac {\mu_i}n).$ This equation is true 
for any graph, regular or not.  \medbreak\noindent{\bf 6.}
Assume now that $G$ is a regular graph of order $n$ and 
degree $d.$ We have $L=dI-A$ and, more generally, 
$L^r=(dI-A)^r=\sum_{i=0}^r(-1)^i{r\choose i}d^{r-i}A^i.$ In general, for any 
square matrix $B,$ we write $B^0=I.$ It follows that 
$tr(L^r)=\sum_{i=0}^r(-1)^i{r\choose i}d^{r-i}tr(A^i).$ 
\medbreak\section*{\Large\bf2.  An exact series formula 
for the log-complexity of a regular graph in terms of closed 
walks} \medbreak\noindent Spanning trees of a graph are 
typically numerous and diverse.  By contrast, walks in a 
graph are just about the easiest to grasp.  Our next 
result expresses the log-complexity of a regular graph as 
an infinite alternating series that involves closed walks.  
Closed walks are traces of the adjacency matrix, and 
while they are particularly intuitive and easy to use, there are other 
meaningful symmetric functions of eigenvalues that can 
be used instead; see [5].  

\medbreak\noindent {\bf Proposition 1} {\em If $G$ is a graph of order $n$ and degree $d$, with $d<n/2$, then the log-complexity of the complement $\bar {G}$ is expressed in terms of $w_k(G
)$, the number of closed walks with $k$ edges in $G$, as follows:}
\medbreak
\centerline{$ln(t(\bar {G}))=$$ln(n^{-2}(n-d)^n)+\sum_{k=2}^{\infty}
(-1)^{k-1}\frac {w_k(G)}{k(n-d)^k}.$}
\medbreak\noindent {\bf Proof} The proof rests on series 
expansions.  We start with the formula for $t(\bar {G})$ in {\bf 5.  }
above, and use all six expressions as needed.  
\medbreak\noindent
$ln(t(\bar {G}))=(n-2)ln(n)+\sum_{i=1}^nln(1-\frac {\mu_i}n)=(n-2
)ln(n)-\sum_{i=1}^n(\sum_{r=1}^{\infty}\frac 1r(\frac {\mu_i}n)^r
)$ 
$=(n-2)ln(n)-(\sum_{r=1}^{\infty}\frac {tr(L^r)}{rn^r}).$ ~ (1) \medbreak\noindent
Observe that, since $0\leq\frac {\mu_i}n<1,$ $\forall i$ the series in (1) 
converges.  Focus on this last series, use the content of 
{\bf 6.},  and change the order of summation.  This yields 
$\sum_{r=1}^{\infty}\frac {tr(L^r)}{rn^r}=\sum_{r=1}^{\infty}\frac 
1{rn^r}(\sum_{k=0}^r(-1)^k{r\choose k}d^{r-k}tr(A^k))=(\sum_{r=1}^{
\infty}\frac 1{rn^r}d^r)tr(A^0)$ 
\vskip.3cm\noindent\ 

$+\sum_{k=1}^{\infty}(\sum_{r=k}^{\infty}\frac 1{rn^r}{r\choose k}
d^{r-k})(-1)^ktr(A^k)=-ln(1-\frac dn)tr(A^0)+\sum_{k=1}^{\infty}\frac {
(-1)^k}{k(n-d)^k}tr(A^k).$ 
\vskip.3cm\noindent
The last sign of equality is explained by making use of 
the identity 
$\sum_{r=k}^{\infty}\frac 1{rn^r}{r\choose k}d^{r-k}=\frac 1{n^k}
(\sum_{s=0}^{\infty}\frac {{{k+s}\choose s}}{k+s}(\frac dn)^s)=n^{
-k}k^{-1}(1-\frac dn)^{-k}=\frac 1{k(n-d)^k}.$ 
Substituting this information into the expression for 
$ln(t(\bar {G}))$ found above, and using {\bf 1.} to introduce the closed 
walks for the traces that arise, we finally obtain 
$ln(t(\bar {G}))=(n-2)ln(n)-(\sum_{r=1}^{\infty}\frac {tr(L^r)}{r
n^r})=(n-2)\cdot ln(n)+n\cdot ln(1-\frac dn)$ 
$+\sum_{k=1}^{\infty}\frac {(-1)^{k-1}tr(A^k)}{k(n-k)^k}=ln(n^{-2}
(n-d)^n)+\sum_{k=2}^{\infty}(-1)^{k-1}\frac {w_k(G)}{k(n-d)^k}.$ 
Since $d$ is the largest eigenvalue of $A$ in absolute value, 
the last series converges when $d<n-d$, which retricts 
$d<n/2$, as enunciated.  This is the expression we 
sought.  \medbreak\noindent It is of some interest to 
assess the speed of convergence of the series in 
Proposition 1.  We first examine an example.  
\medbreak\noindent {\bf Example 1} We take $G$ to be the 
Petersen graph.  Since $G$ is a strongly regular graph 
with $n=10$ and $d=3$, the eigenvalues of the adjacency 
matrix and of the Laplacian are well-known.  We can, 
therefore, directly evaluate the complexity of $\bar {G}$ and 
there is no need for any series expansion.  The point of 
the exercise is two-fold:  we want to check that the 
series expansion gives the correct answer, and we also 
want to examine the speed of convergence of the series.  
\medbreak\noindent The adjacency matrix $A$ of $G$ has 
eigenvalues 1, -2, 3 of respective multiplicities 5, 4, 1.  
The Laplacian of $\bar {G}$ (which is also strongly regular) has 
eigenvalues 8, 5, 0 of multiplicities 5, 4, 1.  It follows 
that $t(\bar {G})=\frac {8^5\cdot 5^4}{10}$.  The number of closed $
k-$walks in $G$ is 
$w_k=5\cdot 1^k+4\cdot (-2)^k+1\cdot 3^k.$ According to Proposition 1, 
$t(\bar {G})=ln(n^{-2}(n-d)^n)+\sum_{k\geq 2}\frac {(-1)^{k-1}w_k}{
k(n-d)^k}.$ Substituting in the 
$w_k$ and using the expansion of the lagarithm series we 
obtain 
\vskip.3cm\noindent
$ln(\bar {G}$$)=ln(10^{-2}\cdot 7^{10})+\sum_{k\geq 2}(-1)^{k-1}\frac {
(5+4(-2)^k+3^k)}{k\cdot 7^k}=$ 
\vskip.3cm\noindent
$ln(10^{-2}\cdot 7^{10})+5(ln(1+\frac 17)-\frac 17)+4(ln(1-\frac 
27)+\frac 27)+ln(1+\frac 37)-\frac 37=$ 
\vskip.3cm\noindent
$ln(10^{-2}\cdot 7^{10})+ln[(1+\frac 17)^5(1-\frac 27)^4(1+\frac 
37)]=ln(\frac {8^5\cdot 5^4}{10}),$ as 
anticipated.  \medbreak\noindent The value $ln(\bar {G})=14.53237$ 
is approximated by the first six partial sum of the 
series given by Proposition 1 as follows:  14.85393, 
14.54781, 14.54781, 14.53219, 14.53362, 14.53221.  As is 
evident from this, on the log-scale an approximation 
obtained by using closed walks of length 4 or less yields 
the correct answer in the first three decimal places.  A 
more detailed look at such approximations takes place in 
the sections that follow.  
\medbreak\noindent {\bf Theorem 1} {\em If $G$ is a graph of order $
n$ and degree $d$, with $d<n/2$, then $t(\bar {G})$ is 
equal to the closest integer to}
\medbreak
\centerline{$n^{-2}(n-d)^nexp(\sum_{k=2}^B\frac {(-1)^{k-1}w_k}{k
(n-d)^k}).$}
\medbreak\noindent
{\em Here $B$ is the smallest natural number with the property 
that }
\medbreak
\centerline{$n^{-2}(n-d)^n\cdot |exp(\sum_{k=2}^m\frac {(-1)^{k-1}
w_k}{k(n-d)^k})-exp(\sum_{k=2}^B\frac {(-1)^{k-1}w_k}{k(n-d)^k})|
<\frac 12,$}
\medbreak\noindent
{\em for all} $m>B.$
\medbreak\noindent
{\bf Proof}  The 
series that appears in Proposition 1 is convergent; it is not always alternating, since some of its terms can be zero.  If we express the series as $\sum_ia_i$ 
then evidently $a_i\rightarrow 0$ as $i\rightarrow\infty .$ But the complexity 
$t(\bar {G})=exp(\sum_ia_i)$ is an {\em integer}.  This tells us that we 
can identify $t(\bar {G})$ by only using the first finite number 
of terms in the series.  Indeed, since the sequence is 
convergent it is also Cauchy and we can stop summing 
when we reach {\em consistent\/} diminishing returns of less 
than $\frac 12$ in the finite product $\prod_{i=1}^mexp(a_i)$, which now 
unambiguously identifies $t(\bar {G}).$ An exact formulation of 
this argument appears in the statement of Theorem 1. 
\medbreak\noindent
{\bf Remarks}
\vskip.3cm\noindent
{\bf 1.} When in the presence of a regular graph $H$ of degree $r$ 
with $n$ vertices it is helpful to always be aware of the 
bounds $L\leq t(H)\leq U$. Here the sharp lower bound $L$ (due to 
McKay [1]) 
takes the value $U=n^{-1}(2r)^kz^{n-1-k}$, with $k$ being the 
integral part of $\frac {n-2}2$ and $z=\frac {nr-2rk}{n-1-k}$; McKay showed that 
this lower bound is {\em always achived\/} by a regular graph 
with $n$ vertices regular of degree $r$. There are many 
choices for the upper bound $U.$ None known to the authors 
are always sharp, but we may use $U=n^{n-2}(\frac r{n-1})^{n-1}$, for 
instance. The integers in the interval $[L,U]$ are therefore the only 
ones that can possibly occur as values for $t(\bar {G})$ in 
Theorem 1. 
\vskip.3cm\noindent
{\bf 2.} In potential applications of Theorem 1 accurate
estimation of the natural number $B$ may become important. 
We believe that this is possible to accomplish by a more 
careful study of the closed walks $w_k(G)$, at least for 
specialized classes of graphs $G.$ The next section 
examines the case of bipartite graphs, for instance, in 
which case bound $B$ (due to monotonicity) is easier to address. We 
shall not deal with the full generality of this issue in 
this paper.  
\medbreak\noindent
\section*{\Large\bf3.  Bounds on complexity} 
\medbreak\noindent There are many upper bounds on graph 
complexity, mostly based on variants of the 
geometric-arithmetic mean inequality and the 
log-concavity of the determinant of a positive definite 
matrix; see $ $[1, 2, 3, 4] and [9].  Lower bounds are rare 
and typically more difficult to obtain; we mentioned 
at the end of the previous section the 
remarkable lower bound found in [1]. We start with 
establishing lower bounds based on the results presented 
in Section 2.  \medbreak\noindent For $p$ a natural number 
and $x$ a vector, we write $|x|_p=(\sum_i|x_i|^p)^{1/p}$ for the $
l_p$ 
norm of $x.$ From inequalities on $l_p$ norms it is known 
and easy to check that $|x|_k\leq |x|_m,$ for $1\leq m\leq k.$ If $
x$ is 
the vector of nonzero eigenvalues of the Laplacian $L$ of 
the connected graph $G$, then it is clear that 
$|x|_k=(tr(L^k))^{1/k}.$ If graph $G$ is of order $n$ and degree $
d$ 
then we may also readily calculate that $trL=nd,$ 
$tr(L^2)=n(d^2+d)$ and $tr(L^3)=nd^3+3nd^2-6\Delta ,$ where $\Delta$ 
stands for the number of triangles in $G.$ We freely use 
these expressions in the remainder of this section.  
\medbreak\noindent As written in Section 2 at the begining 
of the proof of Proposition 1, and by using the $l_p$ norm 
inequalities written above with $m=2,$ we obtain 
\medbreak\noindent
$ln(t(\bar {G}))=(n-2)ln(n)-\sum_{k\geq 1}\frac {tr(L^k)}{kn^k}\geq 
(n-2)ln(n)-\frac {trL}n-\sum_{k\geq 2}\frac {(tr(L^2))^{k/2}}{kn^
k}.$ 
For simplicity let 
$y=\frac {(tr(L^2))^{1/2}}n=\frac {(n(d^2+d))^{1/2}}n=(\frac {d(d
+1)}n)^{1/2}$.  We may now 
write, for $y<1,$ $\sum_{k\geq 2}\frac {(tr(L^2))^{k/2}}{kn^k}=\sum_{
k\geq 2}\frac {y^k}k=-ln(1-y)-y$.  
This yields \medbreak\noindent
$ln(t(\bar {G}))\geq (n-2)ln(n)-d+\sqrt {\frac {d(d+1)}n}+ln(1-\sqrt {\frac {
d(d+1)}n})$, for 
$d(d+1)<n.$ \medbreak\noindent If $\bar {d}$ is the degree of $\bar {
G}$ we 
have $d+\bar {d}=n-1.$ This allows us to express the above 
inequality as 
$ln(t(\bar {G}))\geq (n-2)ln(n)-(n-1-\bar {d})+\sqrt {\frac {(n-1
-\bar {d})(n-\bar {d})}n}+ln(1-$$\sqrt {\frac {(n-1-\bar {d})(n-\bar {
d})}n}),$ 
which is subject to convergence restriction 
$\frac {(n-1-\bar {d})(n-\bar {d})}n<1.$ We summarize as follows, using the 
notation $exp(x)$ to denote the exponential function 
commonly written as $e^x.$ \medbreak\noindent {\bf Proposition 2 }
{\em If $G$ is a graph of order $n$ regular of degree $
d$ satisfying the restriction $(n-1-d)(n-d)<n,$ then $G$ has at least
$n^{n-2}\cdot (1-\sqrt {\frac {(n-1-d)(n-d)}n})\cdot exp(-(n-1-d)
+\sqrt {\frac {(n-1-d)(n-d)}n})$ spanning trees.}  
\medbreak\noindent The degree restriction 
in Proposition 2 is rather severe.  It basically requires 
that the degree $d$ of the graph $G$ be within about a 
square root of $n$ of the degree of the complete graph of 
order $n,$ that is, $d\geq n-\sqrt {n}.$ We show next how this 
restriction can be controlled in large measure by 
expanding the series in powers of $tr(L^m)$ rather than 
simply $tr(L^2).$ This will bring into focuss features of the 
graph other than its degree, such as cycles of higher 
order.  \medbreak\noindent {\bf Theorem 2} {\em If $G$ is a regular 
graph of order $n$ with its Laplacian $L$ satisfying the 
inequality $tr(L^m)<n^m$ for some integer $m\geq 2$, then the 
complexity of the complement $\bar {G}$ verifies the inequality 
\vskip.3cm 
\centerline{$t(\bar {G})\geq n^{n-2}(1-\frac {(tr(L^m))^{\frac 1m}}
n)exp(-\sum_{k=1}^{m-1}\frac {[tr(L^k)-(tr(L^m))^{k/m}]}{kn^k}${\em ).}}
\vskip.3cm\noindent
The inequality becomes equality as $m\rightarrow\infty .$ }
\medbreak\noindent {\bf Proof}  Using the $l_p$ inequalities, for 
$2\leq m\leq k$ we may generally write 
$ln(t(\bar {G}))=(n-2)ln(n)-\sum_{k\geq 1}\frac {tr(L^k)}{kn^k}\geq 
(n-2)ln(n)-\sum_{k=1}^{m-1}\frac {tr(L^k)}{kn^k}-\sum_{k\geq m}\frac {
(tr(L^m))^{k/m}}{kn^k}.$ 
\medbreak\noindent By setting $y=\frac {(tr(L^m))^{1/m}}n$ the above 
inequality may expressed in the form \medbreak\noindent
$ln(t(\bar {G}$$))\geq (n-2)ln(n)+\sum_{k=1}^{m-1}\frac {[(tr(L^m
))^{k/m}-tr(L^k)]}{kn^k}+ln(1-y)$, 
with $0\leq y<1.$ \medbreak\noindent The restriction $0\leq y<1$ 
is equivalent to $tr(L^m)<n^m.$ As written in {\bf 6.}  Section 1, 
with $w_i$ standing for the number of closed $i-$walks, 
$tr(L^m)=\sum_{i=0}^m(-1)^i{m\choose i}d^{m-i}w_i=nd^m+{m\choose 
2}nd^{m-1}-{m\choose 3}w_3d^{m-3}\pm\cdots$ 
which, for sufficiently large fixed $n$ and sufficiently 
small fixed $m$, may conveniently be viewed as a 
polynomial in $d.$ In this asymptotic sense, examining 
just the leading power in $d,$ the inequality $tr(L^m)<n^m$ 
reduces to $d^m<n^{m-1},$ or $d<n^{\frac {m-1}m}.$ It is evident now 
that this last inequality does not actually restrict $d;$ 
for instance, the typical restriction $d<\frac n2$ is verified by 
taking $m$ such that $2^m<n.$ As $m\rightarrow\infty$ we simply 
recapture the incipient content of Proposition 1 as it 
appears in (1) of Section 2.  This ends the proof.  
\medbreak\noindent We now study in further detail the 
case $m=3$ of Theorem 2 since it provides a lower bound 
on complexity in terms of both the degree as well as 
the number of triangles in the graph.  We saw that 
$tr(L^3)=nd^3+3nd^2-6\Delta =nd^2(d+3)-6\Delta ,$ with $\Delta$ signifying 
the number of trangles in the grap $G;$ observe that 
$w_3=\Delta .$ Moreover, simple counting shows that if $d$ 
(respectively $\bar {d}$) and $\Delta$ (respectively $\bar{\Delta}$) denote the degree 
and the number of triangles in $G$ (respectively $\bar {G}$), then 
we have $d+\bar {d}=n-1$ and $\Delta +\bar{\Delta }={n\choose 3}-\frac {
nd\bar {d}}2$; see also [8].  
To simplify notation, write 
$s=n^{-1}\cdot (tr(L^3))^{\frac 13}=n^{-1}\cdot (nd^3+3nd^2-6\Delta 
)^{\frac 13}.$ We deduce 
from Theorem 2 that 
\vskip.3cm
$t(\bar {G}$$)\geq n^{n-2}\cdot (1-s)\cdot exp(s-d+\frac {s^2}2-\frac {
d(d+1)}{2n})$.~ (3) 
\medbreak\noindent Since we are concerned with the graph 
$\bar {G}$, specifically $t(\bar {G}),$ it is helpful to express $
d$ and $s$ 
solely in terms of features of $\bar {G}$ such as $\bar {d}$ and $
\bar{\Delta }.$ We 
have $tr(L^3)=n(n-1-\bar {d})^2(n+2-\bar {d})-6({n\choose 3}-\frac {
n\bar {d}(n-1-\bar {d})}2-\bar{\Delta }).$ In 
summary:  \medbreak\noindent 
{\bf Proposition 3} {\em If $G$ is a graph of 
order $n$ regular of degree $d$ and having $\Delta$ triangles, then 
\vskip.3cm 
\centerline{$t(G)\geq n^{n-2}\cdot (1-s)\cdot exp(s-(n-1-d)+\frac {
s^2}2-\frac {(n-d)(n-d-1)}{2n}),$ }
\vskip.3cm\noindent
where $s$ is defined by 
$n^3s^3=n(n-1-d)^2(n+2-d)-6({n\choose 3}-\frac {nd(n-1-d)}2-\Delta 
).$ The inequality holds true whenever $0\leq s<1.$ }
\medbreak\noindent {\bf Example 2} Consider the graph $H$ with 
10 vertices, labeled 0,1,\ldots ,9 regular of degree 3.  Graph $H$ 
has edges 13, 13, 23, 14, 26, 35, 45, 56, 47, 68, 79, 70, 89, 
80, 90.  We observe that $H$ has 3 triangles.  To start 
with, a direct calculation shows that $G=\bar H$ has 2080524 
spanning trees.  Our interest is in examining the lower 
bound on the complexity of the graph $G=\bar {H}$ as 
highlighted in Proposition 3.  By setting $L_H$ as the 
Laplacian of $H$ we verify that 
$s=\frac {(tr(L_H^3))^{1/3}}n=\frac {\sqrt[3]{522}}{10}=0.8051748
<1,$ which allows the 
application of Proposition 3.  See also (3) for additional 
clarity.  On the log scale we obtain a lower bound of 
14.31436 and may threfore write 
$14.54813=log(t(G))>14.31436.$ Foregoing the log scale, 
the value of the lower bound turns out to be 1646819 
which is indeed less than the true complexity of 
2080524. [We mention in passing that McKay's formula yields a quite low lower bound of 165113 for
the complexity of $G$. This simply means that the graph $G$ has a much higher complexity than the lowest possible complexity graph with 10 vertices and degree 6.]
\vskip.3cm\noindent
It might be interesting to also point out that a lower bound 
for $t(G)$ cannot be obtained by using just the degree, as 
in Proposition 2, since in this example 
$3\cdot 4=d(d+1)<n=10$ does not hold true.  
\medbreak\noindent Presence of triangles in graphs is a 
well-studied problem.  Proposition 3 suggests that graphs 
of maximal complexity among all graphs of a given order 
and specified degree are likely found among those that 
have a minimal number of triangles.  In particular, since 
there is considerable understanding of the structure of 
regular graphs with a minimal number of triangles, cf.  
[11] and [12], this reflects favorably in identifying infinite 
families of graphs of maximal complexity by way of 
Proposition 3 and Proposition 2 above.  Large classes of 
graphs with a minimual number of triangles are 
described in Theorem 1.6 of [12].  A more restricted but 
relatively simple construction appears also in [11].  We 
explain the details.  Let $k$ and $l$ be integers such that 
$k>l\geq 0.$ Start with a complete bipartite graph 
$K_{2k+l,2k+l}$ with vertex set $\{x_1,\ldots ,x_{2k+1}\}$ and 
$\{y_1,\ldots ,y_{2k+l}\}$.  Remove a $(l+1)-$factor from the graph 
induced by set $x_1,\ldots ,x_k,y_1,\ldots ,y_k$ and an $l$-factor from 
the graph induced by $x_{k+1},\ldots ,x_{2k+l},y_{k+1},\ldots ,y_{
2k+l}$.  Join 
$x_1,\ldots ,x_k,y_1,\ldots ,y_k$ to a new vertex $z.$ Denote by $
g(k,l)$ the 
family of graphs so obtained.  An element of $g(k,l)$ is a 
regular graph of degree $2k$ with $4k+2l+1$ vertices.  It is 
shown in [11, Theorem 2.1] that for $k\geq 2^{20}$ and 
$k\geq 2l+6\sqrt {10l}+1$ a graph in $g(k,l)$ is the sole graph with 
$4k+2l+1$ vertices of degree $2k$ having a minimal number 
of triangles (exactly $k(k-l-1)$ triangles) among all graphs 
with the same number of vertices and of the same 
degree.  When viewed in the context of Proposition 1 and 
Proposition 2 above, the results contained in Theorem 1.6 
of [12] and Theorem 2.1 of [11], provide us with infinite 
families of graphs that have few short closed walks and 
would therefore also have high complexity.  As explained 
in the Introducion, this is the desirable feature that we 
want in facilitating neuronal signal transmission.  
\medbreak\noindent\section*{\Large\bf4.  Complements of 
bipartite graphs} \medbreak\noindent As is well-known, a 
bipartite graph has no closed walks of odd length, and is 
characterized by this property.  We use the results in 
the previous two sections to investigate the complexity 
of graphs that are complements of bipartite graphs.  Let 
$G$ be a bipartite graph of order $n,$ regular of degree $d.$ 
We remind that $w_k(G)$ denotes the number of closed 
$k-$walks in $G.$ When the presence of $G$ is understood we 
simply write $w_k$ for $w_k(G).$ As mentioned, the bipartite 
assumption on $G$ forces $w_k(G)=0$ for all odd $k\geq 1.$ 
\medbreak\noindent {\bf Theorem 3} {\em If $G$ is a bipartite graph of order 
$n$ regular of degree $d$, and $m,$ $k$ are positive integers,
then $a(n,d,m)\leq t(\bar {G}$$)\leq b(n,d,k),$ where 
\vskip.3cm\noindent
$a(n,d,m)=(n-d)^n\cdot n^{-2}\cdot\sqrt {1-y^2}\cdot exp(-\sum_{1
\leq s<m}\frac {(w_{2s}-(w_{2m})^{s/m})}{2s(n-d)^{2s}}),$ 
\vskip.3cm\noindent
$b(n,d,k)=(n-d)^n\cdot n^{-2}\cdot exp(-\sum_{s=1}^k\frac {w_{2s}}{
2s(n-d)^{2s}})$ and 
$y=\frac {(w_{2m})^{1/2m}}{n-d}.$ The lower bound holds true whenever 
$y<1.$ When $m\rightarrow\infty$ or when $k\rightarrow
\infty$ the respective 
inequalities become equalities.  \medbreak}\noindent {\bf Proof }
For such $G$ Proposition 1 takes the form 
\vskip.3cm\noindent
\centerline{$ln(t(\bar {G}))=$$ln(n^{-2}(n-d)^n)-\sum_{s=1}^{\infty}\frac {
w_{2s}(G)}{2s(n-d)^{2s}}.$}
\vskip.3cm\noindent
From this, the choice of $b(n,d,k)$ immediately follows.  
We now explain how the lower bound $a(n,d,m)$ is 
achieved.  Relying on {\bf 1.}  and {\bf 2.}  in Section 1, 
$w_{2s}:=w_{2s}(G)=tr(A^{2s}),$ where $A$ is the adjacency matrix 
of $G$.  The eigenvalues of $A^2$ are nonnegative since they 
are the squares of the (real) eigenvalues of $A.$ Making 
use of the $l_p$ inequalities we may write 
$tr(A^{2s})\leq (tr(A^{2m}))^{2s/2m},$ for $s\geq m\geq 1.$ With 
$y=\frac {tr(A^{2m})^{1/2m}}{n-d}=\frac {(w_{2m})^{1/2m}}{n-d},$ this yields 
\vskip.3cm\noindent
$\sum_{s=1}^{\infty}\frac {w_{2s}}{2s(n-d)^{2s}}\leq\sum_{1\leq s
<m}\frac {tr(A^{2s})}{2s(n-d)^{2s}}+\sum_{s\geq m}\frac {(tr(A^{2
m}))^{2s/2m}}{2s(n-d)^{2s}}=$ 
\vskip.3cm\noindent
$\sum_{1\leq s<m}\frac {w_{2s}}{2s(n-d)^{2s}}+\sum_{s\geq m}\frac {
y^{2s}}{2s}=\sum_{1\leq s<m}\frac {w_{2s}}{2s(n-d)^{2s}}$ 
\vskip.3cm\noindent
$-\frac 12[ln(1-y)+ln(1+y)]-\sum_{1\leq s<m}\frac {y^{2s}}{2s}$.  \medbreak\noindent
Exponentiating both sides of the inequality yields 
\vskip.5cm\noindent

$t(\bar {G}$$)\geq (n-d)^n\cdot n^{-2}\cdot\sqrt {1-y^2}\cdot exp
(-\sum_{1\leq s<m}\frac {(w_{2s}-(w_{2m})^{s/m})}{2s(n-d)^{2s}})=
a(n,d,m)$ 
\vskip.3cm\noindent
as enunciated.  From the formula in Proposition 1 it 
follows that when $m\rightarrow\infty$ or when $k\rightarrow\infty$ the inequalities 
become equalities.  This ends the proof.  
\medbreak\noindent We illustrate the content of Theorem 3 
by an example.  \medbreak\noindent {\bf Example 3} Consider the 
bipartite graph $G$ on vertices 0,1,\ldots ,9 regular of degree 3 
with parts 1,2,3,4,5 and 6,7,8,9,0.  Edges of $G$ are 17 18 19 
28 29 20 36 39 30 40 46 47 56 57 58.  Direct 
computation shows $t(\bar {G})=2034010.$ For graph $G$ we have 
$w_2=30,$ $w_4=190,$ $w_6=1530,$ \ldots\ We examine the bounds 
for values $(m,k)\in \{(2,2),$ (3,3), (4,4), $(5,5)$, $(6,6)\}$ The 
corresponding values for $(a(n,d,m),b(n,d,k))$ are as 
follows:  (2029504, 2039113), (2033738, 2034698), (2033985, 
2034111), (2034007, 2034025), (2034010, 2034012).  We 
observe that for $m=6$ the lower bound yields the $exact$ 
answer.  It turns out that at $k=7$ the upper bound also 
equals the exact answer.  
\vskip1cm 
\centerline{{\bf Acknowledgement}}
\vskip.3cm\noindent\ 
We are grateful to the National Science Foundation for 
sponsoring this work under the Funding Opportunity NSF 
24-513 Emerging Mathematics in Biology.  
\vskip1cm 

\centerline{{\bf REFERENCES}}
\begin{enumerate}
\item McKay, B.  D.  Spanning trees in regular graphs, 
{\em Eur.  J.  Comb.}, 4, 149-160 (1983) 

\item Das, K.  A.  A sharp upper bound for the number of 
spanning trees of a graph, {\em Graphs Comb.,\/}  23, 625-632 
(2007) 

\item Li, J.,  Shiu, W.C.,  Chang, A.  The number of 
spanning trees of a graph, $Appl.$ $Math.$ $Lett.,$ 23, 286-290 
(2010) 

\item Chung, F.,  Yau, S-T.  Coverings, heat kernels and 
spanning trees, {\em Electron.  J.  Comb.,\/}  6, R12 (1999) 

\item MacDonald, I.  G.  {\em Symmetric functions and Hall polynomials}, Oxford University press, 2015 

\item Brouwer, A.  E.  and Haemers, W.  H.  {\em Spectra of }
{\em graphs,\/} Springer, New York, 2012 

\item van Dam, E.  R.  Graphs with few eigenvalues, {\em PhD dissertation}, Tilburg University, 1996 

\item Radhakrishnan, N.  and Vijayakumar A.  (1994), 
About triangles in a graph and its complement, {\em Discrete mathematics,\/} 131, 205-210 

\item Alon, N.  The number of spanning trees in regular 
graphs, {\em Random structures and algorithms,\/} vol 1 (2), 
175-191 (1990) 

\item Constantine, G.  M.  {\em Combinatorial theory and statistical design}, Wiley, New York, 1987 

\item Lo, A.  S.  L.  (2009) Triangles in regular graphs 
with density below one half, {\em Combinatorics, Probability and Computing,\/} 18, 435-440 

\item Liu, H.,  Pikhurko, O.,  Staden, K.  (2020) The exact 
minimum of triangles in graphs with given order and 
size, {\em Forum of Mathematics, Pi,\/} Vol.  8, e8, 144 pages 
doi:  10.1017/fmp.2020.7 

\item Bohnen, N.,  Prabesh, K.,  Koeppe R.,  Catasus, C.,  
Frey, K.,  Scott, P.,  Constantine, G.,  Albin, R, M\"uller, 
M.  (2021) Regional cerebral cholinergic nerve terminal 
integrity and cardinal motor features in Parkinson's 
disease, {\em Brain communications,\/} vol 3, issue 2, fcab109 

\item Bear, M.,  Connors, B.,  Paradiso, M.  {\em Neuroscience:  Exploring the brain,\/} Fourth edition, Jones and Bartlett 
Learning, Burlington, MA, 2016 

\end{enumerate}

\end{document}